\documentclass[12pt,a4paper]{amsart}

\newtheorem{Result}{Result}[section]
\newtheorem{Thm}  [Result]{Theorem}
\newtheorem{Prop} [Result]{Proposition}
\newtheorem{Lemma}[Result]{Lemma}
\newtheorem{Cor}  [Result]{Corollary}

\newtheorem*{Def}{Definition}
\newtheorem*{Remark}{Remark}
\expandafter\def\expandafter\Remark\expandafter{\Remark\upshape}

\newtheorem*{Question}{Question}

\def\qedbox{\hbox{\vrule\vbox{\hrule width 6pt\vskip6pt\hrule}\vrule}}
\def\qed{\ifvmode\leavevmode\fi
  \unskip\nobreak\hfill\penalty50 \quad \null\nobreak\hfill
  \qedbox{\parfillskip0pt \finalhyphendemerits0 \par}}
\newenvironment{Proof}{\ifdim\lastskip=0pt \medskip\fi
  \noindent{\it Proof.~}\ignorespaces}{\qed\medskip}

\def\infdim{in\-fi\-nite-dimen\-sional}
\def\infcodim{in\-fi\-nite-codi\-men\-sional}

\def\rnbox#1{(\romannumeral#1)}
\def\itemn#1{\item[\rnbox{#1}]}

\def\implies{\DOTSB\;\Rightarrow\;}

\def\N{{\mathbb{N}}}
\let\eps\varepsilon

\let\<\langle
\let\>\rangle

\def\lop{{\mathcal{L}}}

\def\JL{{\mathchoice{\hbox{\it JL\/}}{\hbox{\it JL\/}}%
  {\hbox{\the\scriptfont1 JL\/}}{\hbox{\the\scriptscriptfont1 JL\/}}}}

\let\maparrow\longrightarrow
\def\map#1#2#3{#1\colon#2\,{\maparrow}\,#3}

\def\eqalign#1{\null\,\vcenter{\openup\jot\mathsurround0pt\relax
  \ialign{\strut\hfil$\displaystyle{##}$&$\displaystyle{{}##}$\hfil
      \crcr#1\crcr}}\,}

\begin{document}

\title[The three-space property for subprojective Banach spaces]{On the
  three-space property for subprojective and superprojective Banach~spaces}

\author{Manuel Gonz\'alez}
\address{Facultad de Ciencias, Universidad de Cantabria, 39071~Santander, Spain}

\author{Javier Pello}
\address{Escuela Superior de Ciencias Experimentales y Tecnolog\'\i a,
  Universidad Rey Juan Carlos, 28933~M\'ostoles, Spain}

\begin{abstract}
We introduce the notion of subprojective and superprojective operators
and we use them to prove a variation of the three-space property for
subprojective and superprojective spaces. As an application, we show
that some spaces considered by Johnson and Lindenstrauss are both
subprojective and superprojective.
\end{abstract}

\maketitle

\begin{section}{Introduction}

A Banach space~$X$ is called \emph{subprojective} if every closed
\infdim{} subspace of~$X$ contains an \infdim{}
subspace complemented in~$X$, and $X$ is called \emph{superprojective}
if every closed \infcodim{} subspace of~$X$ is contained in
an \infcodim{} subspace complemented in~$X$; note that
finite-dimensional spaces are trivially both subprojective and
superprojective.
These two classes of Banach spaces were introduced by Whitley~\cite{whitley}
in order to find conditions for the conjugate of an operator to be strictly
singular or strictly cosingular.
More recently, they have been used to obtain some positive solutions to the
perturbation classes problem for semi-Fredholm operators. This problem has a
negative solution in general~\cite{perturbation}, but there are some positive
answers when one of the spaces is subprojective or superprojective
\cite{perturbation-sub-sup} \cite{G-P-Salas}.

Subprojectivity passes on to subspaces and superprojectivity passes on to
quotients, and both are stable under direct sums \cite{oikhberg-spinu}
\cite{superprojective}, but neither of them  is a three-space property:
for every $1 < p < \infty$ there exists a non-subprojective,
non-superprojective space~$X$ with a subspace $M\subseteq X$ such that
$M\simeq X/M\simeq\ell_p$, which is both subprojective and superprojective
\cite[Proposition 2.8]{oikhberg-spinu}
\cite[Proposition 3.2]{superprojective}.
However, slightly stronger hypotheses on $M$ and~$X/M$ do imply the
subprojectivity or superprojectivity of~$X$, as seen below.
For the subprojectivity, it is enough that $X/M$ is subprojective and
that every closed \infdim{} subspace of~$M$ contains an
\infdim{} subspace complemented in~$X$ (which is stronger than
being complemented in~$M$, as the definition of subprojectivity of~$M$ would
require), and there is an equivalent condition for superprojectivity
(see Theorem~\ref{3sp}).
We introduce two classes of operators, namely subprojective and
superprojective operators, which allow to show that a Banach space~$X$
is subprojective or superprojective when certain conditions such as these
are met by a closed subspace~$M$ of~$X$ and its induced quotient~$X/M$.

As an application, these sufficient conditions will be used to prove
that two examples of Banach spaces introduced by Johnson and Lindenstrauss
to study the properties of weakly compactly generated spaces are both
subprojective and superprojective.

We will use standard notation. $X$, $Y$ and~$Z$ will be Banach spaces.
Given a closed subpsace~$M$ of~$X$, we will denote the inclusion of~$M$
into~$X$ by~$J_M$, and $Q_M$ will be the quotient map of~$X$ onto~$X/M$.
A (bounded, linear) operator $T\in\lop(X,Y)$ is said to be strictly singular
if there is no closed \infdim{} subspace~$M$ of~$X$ such that
the restriction $T J_M$ is an isomorphism; it is said to be strictly
cosingular if there is no closed \infcodim{} subspace~$M$ of~$Y$
such that $Q_M T$ is surjective.

\end{section}

\begin{section}{Subprojective and superprojective operators}

\begin{Def}
An operator $T\in\lop(X,Y)$ is subprojective if every closed
\infdim{} subspace~$M$ of~$X$ such that $T J_M$ is an isomorphism
contains a closed \infdim{} subspace~$N$ such that $T(N)$ is
complemented in~$Y$.

An operator $T\in\lop(X,Y)$ is superprojective if every closed
\infcodim{} subspace~$M$ of~$Y$ such that $Q_M T$ is surjective
is contained in a closed \infcodim{} subspace~$N$ such that
$T^{-1}(N)$ is complemented in~$X$.
\end{Def}

Note that a Banach space~$X$ is subprojective (resp., superprojective)
if and only if the identity~$I_X$ is subprojective (resp., superprojective).
Also, strictly singular operators are trivially subprojective, and strictly
cosingular operators are trivially superprojective.

The following result, which is a consequence of
\cite[Lemma~2.2]{improjective}, will be useful at several places.

\begin{Lemma}\label{useful}
Let $X$ and~$Y$ be Banach spaces and let $T\in\lop(X,Y)$ be an operator.
\begin{enumerate}
\itemn1 If $M$ is a closed subspace of~$X$ such that $T J_M$ is an isomorphism
and $T(M)$ is complemented in~$Y$, then $M$ is complemented in~$X$.
\itemn2 If $N$ is a closed subspace of~$Y$ such that $Q_N T$ is surjective
and $T^{-1}(N)$ is complemented in~$X$, then $N$ is complemented in~$Y$.
\end{enumerate}
\end{Lemma}

\begin{Proof}
(i)~Let $N(T)$ be the kernel of~$T$. If $M \cap N(T) = 0$ and $N$ is a closed
subspace of~$Y$ such that $Y = T(M) \oplus N$, then $X = M \oplus T^{-1}(N)$.

(ii)~Let $R(T)$ be the range of~$T$. If $R(T) + N = Y$ and $M$ is a closed
subspace of~$X$ such that $X = M \oplus T^{-1}(N)$, then $T(M)$ is closed
and $Y = T(M) \oplus N$.
\end{Proof}

Subprojective (resp., superprojective) operators are stable under left
(resp., right) composition.

\begin{Prop}\label{stable}
Let $X$, $Y$ and~$Z$ be Banach spaces and let $T\in\lop(X,Y)$ and
$S\in\lop(Y,Z)$ be operators.
\begin{enumerate}
\itemn1 If $S$ is subprojective, then $ST$ is subprojective.
\itemn2 If $T$ is superprojective, then $ST$ is superprojective.
\end{enumerate}
\end{Prop}

\begin{Proof}
(i)~Let $M$ be a closed \infdim{} subspace of~$X$ such that
$S T J_M$ is an isomorphism. Then $T(M)$ is a closed \infdim{}
subspace of~$Y$ and $S J_{T(M)}$ is an isomorphism. Since $S$ is
subprojective, there exists a closed \infdim{} subspace~$N$
of~$T(M)$ such that $S(N)$ is complemented in~$Z$, and then
$(T J_M)^{-1}(N)$ is a closed \infdim{} subspace of~$M$
whose image $(ST) \bigl( (T J_M)^{-1}(N) \bigr) = S(N)$ is complemented
in~$Z$.

(ii)~Let $M$ be a closed \infcodim{} subspace of~$Z$ such that
$Q_M S T$ is surjective. Then $S^{-1}(M)$ is a closed \infcodim{}
subspace of~$Y$ and $Q_{S^{-1}(M)} T$ is surjective. Since $T$ is
superprojective, there exists a closed \infcodim{} subspace~$N$
of~$Y$ containing~$S^{-1}(M)$ such that $T^{-1}(N)$ is complemented in~$X$,
and then $S(N)$ is a closed \infcodim{} subspace of~$Z$
containing~$M$ where $(ST)^{-1}(S(N)) = T^{-1}(N)$ is complemented in~$X$.
\end{Proof}

Applying this result to the identity of a subprojective or superprojective
space yields the following.

\begin{Cor}\label{space-operator}
Let $X$ and~$Y$ be Banach spaces and let $T\in\lop(X,Y)$ be an operator.
\begin{enumerate}
\itemn1 If $Y$ is subprojective, then $T$ is subprojective.
\itemn2 If $X$ is superprojective, then $T$ is superprojective.
\end{enumerate}
\end{Cor}

\begin{Question}
Let $S$, $T\in\lop(X,Y)$.
\begin{enumerate}
\item If $S$ and $T$ are subprojective, is $S+T$ subprojective?
\item If $S$ and $T$ are superprojective, is $S+T$ superprojective?
\end{enumerate}
\end{Question}

The subprojectivity of an embedding and the superprojectivity of a quotient
map appear often enough that it is worth noting the following
characterisation for them.

\begin{Prop}\label{wrt}
Let $X$ be a Banach space and let $Z$ be a closed subspace of~$X$.
\begin{enumerate}
\itemn1 $J_Z$ is subprojective if and only if every closed
\infdim{} subspace of~$Z$ contains an \infdim{}
subspace complemented in~$X$, in which case $Z$ is subprojective.
\itemn2 $Q_Z$ is superprojective if and only if every closed
\infcodim{} subspace of~$X$ containing~$Z$ is contained in an
\infcodim{} subspace complemented in~$X$, in which case
$X/Z$ is superprojective.
\end{enumerate}
\end{Prop}

\begin{Proof}
(i) This is a direct consequence of the definition of subprojective
operator and Lemma~\ref{useful}(i).

(ii) Assume that $Q_Z$ is superprojective and let $M$ be a closed
\infcodim{} subspace of~$X$ containing~$Z$. Then $Q_Z(M)$ is a
closed \infcodim{} subspace of~$X/Z$ and, by hypothesis,
there exists an \infcodim{} subspace~$N$ of~$X/Z$
containing~$Q_Z(M)$ and such that $Q_Z^{-1}(N)$ is complemented in~$X$,
where $Q_Z^{-1}(N)$ contains~$M$ and is still \infcodim.

Conversely, let $M$ be a closed \infcodim{} subspace of~$X/Z$.
Then $Q_Z^{-1}(M)$ is a closed \infcodim{} subspace of~$X$
containing~$Z$, so there exists an \infcodim{} subspace~$N$
containing~$Q_Z^{-1}(M)$ and complemented in~$X$, and then $Q_Z(N)$ is
\infcodim{} in~$X/Z$ and contains~$M$, and $Q_Z^{-1}(Q_Z(N)) = N$
is complemented in~$X$.

Finally, let $M$ be a closed \infcodim{} subspace~$M$ of~$X/Z$.
Then $Q_Z^{-1}(M)$ is contained in an \infcodim{} subspace~$N$
complemented in~$X$, so $M$ is contained in~$Q_Z(N)$, which is closed,
hence complemented in~$X/Z$ by Lemma~\ref{useful}(ii).
\end{Proof}

The following lemmata will be used in Theorem~\ref{3sp} to handle subspaces of
a space depending on its relative position with respect to another subspace.

\begin{Lemma}\label{sub-3sp-lemma}
Let $X$ be a Banach space, let $M$ and $N$ be closed 
subspaces of~$X$ such that $M \cap N = 0$ and $M + N$ is not closed.
Then there exists an automorphism $\map UXX$ such that $U(M) \cap N$ is
\infdim.
\end{Lemma}

\begin{Proof}
Take normalised sequences $(x_n)_{n\in\N}$ in~$M$ and $(y_n)_{n\in\N}$ in~$N$
such that $\|x_n - y_n\| < 2^{-n}$ for every $n\in\N$. Since any weak
cluster point of $(x_n)_{n\in\N}$ must be in $M\cap N = 0$, by passing
to a subsequence \cite[Theorem~1.5.6]{kalton-albiac} we can assume that
$(x_n)_{n\in\N}$ is a basic sequence and that there exists a sequence
$(x^*_n)_{n\in\N}$ in~$X^*$ such that $\<x^*_i,x_j\> = \delta_{ij}$ for
every $i$, $j\in\N$ and $\sum_{n=1}^\infty \|x^*_n\| \, \|x_n - y_n\| < 1$.
Then $K(x) = \sum_{n=1}^\infty \<x^*_n,x\> (x_n - y_n)$ defines an operator
$\map KXX$ with $\|K\| < 1$ and $U = I - K$ is an automorphism on~$X$
that maps $U(x_n) = y_n$ for every $n\in\N$, so $U(M) \cap N$ is \infdim.
\end{Proof}

\begin{Lemma}\label{sup-3sp-lemma}
Let $X$ be a Banach space, let $M$ and $N$ be closed 
subspaces of~$X$ such that $M + N$ is dense in~$X$ but not closed.
Then there exists an automorphism $\map UXX$ such that
$\overline{U^{-1}(M) + N}$ is \infcodim{} in~$X$.
\end{Lemma}

\begin{Proof}
$M + N$ is dense in~$X$ but not closed, so $M^\bot \cap N^\bot = 0$ and
$M^\bot + N^\bot$ is not closed either \cite[Theorem~IV.4.8]{kato}. Take
a normalised sequence $(x^*_n)_{n\in\N}$ in~$M^\bot$ and another sequence
$(y^*_n)_{n\in\N}$ in~$N^\bot$ such that $\|x^*_n - y^*_n\| < 2^{-n}$ for
every $n\in\N$. Since any weak$^*$ cluster point of $(x^*_n)_{n\in\N}$
must be in $M^\bot \cap N^\bot = 0$, by passing to a subsequence
\cite[Lemma 3.1.19]{tauberian} we can assume that $(x^*_n)_{n\in\N}$ is a
basic sequence and find a sequence $(x_n)_{n\in\N}$ in~$X$ such that
$\<x^*_i,x_j\> = \delta_{ij}$ for every $i$, $j\in\N$ and
$\sum_{n=1}^\infty \|x_n\| \, \|x^*_n - y^*_n\| < 1$.
Then $K(x) = \sum_{n=1}^\infty \< x^*_n - y^*_n, x \> x_n$ defines an operator
$\map KXX$ with $\|K\| < 1$ and $U = I - K$ is an automorphism on~$X$
whose conjugate maps $U(x^*_n) = y^*_n$ for every $n\in\N$, so
$U^*(M^\bot) \cap N^\bot = (U^{-1}(M) + N)^\bot$ is \infdim{}
and $\overline{U^{-1}(M) + N}$ is \infcodim.
\end{Proof}

\begin{Thm}\label{3sp}
Let $X$ be a Banach space and let $Z$ be a closed subspace of~$X$.
\begin{enumerate}
\itemn1 If $J_Z$ and $Q_Z$ are both subprojective, then $X$ is subprojective.
\itemn2 If $J_Z$ and $Q_Z$ are both superprojective, then $X$ is
superprojective.
\end{enumerate}
\end{Thm}

\begin{Proof}
(i)~Let $M$ be a closed \infdim{} subspace of~$X$.
If $M\cap Z$ is \infdim, then it contains another
\infdim{} subspace complemented in~$X$ by the hypothesis on~$J_Z$
and Proposition~\ref{wrt}(i).

Otherwise, if $M\cap Z$ is finite-dimensional, we can assume that
$M\cap Z = 0$ by passing to a further subspace if necessary.
If $M + Z$ is closed, then $Q_Z J_M$ is an isomorphism and $M$ contains an
\infdim{} subspace complemented in~$X$ by the hypothesis on~$Q_Z$.

We are left with the case where $M\cap Z = 0$ and $M + Z$ is not closed.
By Lemma~\ref{sub-3sp-lemma}, there exists an automorphism $\map UXX$ such
that $U(M) \cap Z$ is \infdim. Let $N$ be an
\infdim{} subspace of $U(M) \cap Z$ complemented in~$X$, again
by Proposition~\ref{wrt}(i); then $U^{-1}(N) \subseteq M$ and is still
complemented in~$X$.

(ii)~Let $M$ be a closed \infcodim{} subspace of~$X$.
If $\overline{M + Z}$ is \infcodim, then it is contained in
another \infcodim{} subspace complemented in~$X$ by the
hypothesis on~$Q_Z$ and Proposition~\ref{wrt}(ii).

Otherwise, if $\overline{M + Z}$ is finite-codimensional, we can assume that
$\overline{M + Z} = X$ by enlarging~$M$ with a finite-dimensional subspace
if necessary. If $M + Z$ is closed, so $M + Z = X$, then $Q_M J_Z$ is
surjective and $M$ is contained in an \infcodim{} subspace
complemented in~$X$ by the hypothesis on~$J_Z$.

We are left with the case where $M + Z$ is dense in~$X$ but not closed.
By Lemma~\ref{sup-3sp-lemma}, there exists an automorphism $\map UXX$ such
that $\overline{U^{-1}(M) + Z}$ is \infcodim{} in~$X$.
Let $N$ be an \infcodim{} subspace complemented in~$X$
such that $\overline{U^{-1}(M) + Z} \subseteq N$, again by
Proposition~\ref{wrt}(ii); then $M \subseteq U(N)$, which is still
\infcodim{} and complemented in~$X$.
\end{Proof}

Theorem~\ref{3sp} implies a variation of the 3-space property for
subprojectivity and superprojectivity. Given a Banach space~$X$ and a
closed subspace~$Z$ of~$X$, the inclusion~$J_Z$ is subprojective if and
only if every closed \infdim{} subspace of~$Z$ contains an
\infdim{} subspace complemented in~$X$ by Proposition~\ref{wrt}(i),
and this is stronger than being subprojective; on the other hand, for $Q_Z$
to be subprojective, it is sufficient (but not necessary) that $X/Z$ be
subprojective, by Corollary~\ref{space-operator}. So, if $X/Z$ is
subprojective and 
$J_Z$ is subprojective, then $X$ is subprojective. Similarly, if $Z$ is
superprojective and $Q_Z$ is superprojective, then $X$ is superprojective.

Also, Theorem~\ref{3sp} is not a characterisation. While $J_Z$ must be
subprojective if $X$ is subprojective by Corollary~\ref{space-operator},
$Q_Z$ need not be. For instance, take any surjection
$\map{T}{\ell_1}{C([0,1]})$ and define the operator
$\map{Q}{\ell_1\oplus\ell_1}{C([0,1])}$ as $Q(x,y) = x + Ty$. Then $Q$ is
also clearly surjective and an isomorphism on its first component, but
$\ell_1$ cannot contain any subspace whose image by~$Q$ is complemented
in~$C([0,1])$.

Related to Theorem~\ref{3sp}, as mentioned in the introduction, there exists a
non-subprojective, non-superprojective space~$X$ with a subspace
$M\subseteq X$ such that $M \simeq X/M \simeq \ell_p$, which is both
subprojective and superprojective \cite[Proposition~2.8]{oikhberg-spinu}
\cite[Proposition~3.2]{superprojective}, where $Q_M$ is strictly singular
\cite[Theorem~6.4]{kalton-peck}, hence subprojective, while $X$ is not,
and $J_M$ is strictly cosingular \cite[Theorem~6.4]{kalton-peck}, hence
superprojective, while $X$ is not. Thus, $J_M$ cannot be subprojective
although $M$ is, and $Q_M$ cannot be superprojective although $X/M$ is.
$$\eqalign{
  \hbox{$J_Z$, $Q_Z$ subprojective}
    & \implies \hbox{$X$ subprojective} \implies \cr
    & \implies \hbox{$J_Z$ subprojective}
      \implies \hbox{$Z$ subprojective} \cr
}$$
A similar situation holds for superprojectivity.
$$\eqalign{
  \hbox{$J_Z$, $Q_Z$ superprojective}
    & \implies \hbox{$X$ superprojective} \implies \cr
    & \hskip-5em \implies \hbox{$Q_Z$ superprojective}
      \implies \hbox{$X/Z$ superprojective} \cr
}$$

As a particular case of Theorem~\ref{3sp}, we have the following.

\begin{Cor}\label{main-cor}
Let $X$ be a Banach space such that
\begin{enumerate}
\itemn1 every closed \infdim{} subspace of~$X$ contains an
\infdim{} subspace complemented in~$X^{**}$; and
\itemn2 $X^{**}/X$ is subprojective.
\end{enumerate}
Then $X^{**}$ is subprojective.
\end{Cor}

In~\cite{j-sum-subprojective}, it is proved that, under certain conditions,
the $J$-sum $J(\Phi)$ of Banach spaces, as defined by Bellenot~\cite{Jsum},
and its bidual $J(\Phi)^{**}$ are subprojective. As a particular case, for
every separable subprojective space~$X$ there exists a separable Banach
space~$J(\Phi)$ such that $J(\Phi)^{**}$ is subprojective and
$J(\Phi)^{**}/J(\Phi)$ is isomorphic to~$X$. An essential part of the
proof is to show that any closed \infdim{} subspace of~$J(\Phi)$
contains a further \infdim{} subspace complemented
in~$J(\Phi)^{**}$, which is to say that the inclusion of $J(\Phi)$
in~$J(\Phi)^{**}$ is subprojective \cite[Theorem~5.2(i)]{j-sum-subprojective}.
That result is, in fact, a particular case of Corollary~\ref{main-cor}.

A related result was proved by Argyros and
Raikoftsalis~\cite{argyros-raikoftsalis}.
Let $Y$ be a separable reflexive space. Then:
\begin{itemize}
\item For every $1 < p < \infty$, there is a separable reflexive space
$X_p(Y)$ that is hereditarily complemented~$\ell_p$, hence subprojective,
and $Y$ is a quotient of~$X_p(Y)$.
\item $Y$ is a quotient of a separable hereditarily-$c_0$ space $X_0(Y)$.
\end{itemize}
Of course, in both cases the kernel of the quotient map is subprojective.

\begin{Remark}
In Corollary~\ref{main-cor}, condition~(i) can be replaced by
\begin{enumerate}
\item[(i')] \textit{every closed \infdim{} subspace of~$X$
contains an \infdim{} reflexive subspace complemented in~$X$.}
\end{enumerate}
\end{Remark}

\end{section}

\begin{section}
  {Subprojectivity and superprojectivity of~the~Johnson-Lindenstrauss space}

Here we apply the results in the previous section to study the Banach spaces
introduced by Johnson and Lindenstrauss in~\cite[Examples 1 and~2]{wcg}.

Let $\Gamma$ be a set with the cardinality of the continuum and let
$\{\, N_\gamma : \gamma\in\Gamma \,\}$ be a family of infinite subsets
of~$\N$ such that $N_\gamma \cap N_{\gamma'}$ is finite if $\gamma\ne\gamma'$.
For each $\gamma\in\Gamma$, let $\phi_\gamma \in \ell_\infty$ be the
characteristic function of~$N_\gamma$.

Let $V = c_0 \subset \ell_\infty$ and let $\JL_0$ be the linear span of
$V \cup \{\, \phi_\gamma : \gamma\in\Gamma \,\}$ in~$\ell_\infty$
endowed with the norm
$$\biggl\| y + \sum_{i=1}^k a_i \phi_{\gamma_i} \biggr\|_\JL =
  \max \biggl\{\,
      \biggl\| y + \sum_{i=1}^k a_i \phi_{\gamma_i} \biggr\|_\infty,
      \bigl\| (a_i)_{i=1}^k \bigr\|_2 \,\biggr\}$$
for every $y\in V$ and $(\gamma_i)_{i=1}^k$ with $\gamma_i\ne\gamma_j$
if $i\ne j$.

The space~$\JL$ is defined as the completion of $(\JL_0, \|\cdot\|_\JL)$.
These are some of its properties.

\begin{Thm}\label{wcg}\cite[Example~1]{wcg}
\begin{enumerate}
\itemn1 $V$ is a subspace of~$\JL$ isometric to~$c_0$ and $\JL/V$ is
isometric to~$\ell_2(\Gamma)$.
\itemn2 Weakly compact subsets of~$\JL$ are separable. Hence every
reflexive subspace of~$\JL$ is separable and $V$ is not complemented.
\itemn3 $\JL^*$ is isomorphic to $\ell_1 \oplus \ell_2(\Gamma)$.
\end{enumerate}
\end{Thm}

A Banach space~$X$ is weakly compactly generated (WCG, in short) if there
exists a weakly compact subset of~$X$ that generates a subspace that is
dense in~$X$. Clearly, separable spaces and reflexive spaces are WCG,
but $\JL$ is not by property~(ii), as it is not separable.
Hence being WCG is not a three-space property~\cite{wcg}.

Note that, while we treat~$\JL$ as a unique space, there are different
$\JL$ spaces depending on the choice of the family
$\{\, N_\gamma : \gamma\in\Gamma \,\}$, and that the resulting spaces
may not be isomorphic~\cite{jl-weak-star}.
However, the properties in Theorem~\ref{wcg}, and those proved below, are
common to all possible $\JL$ spaces obtained this way.

\begin{Cor}
Every \infdim{} reflexive subspace of~$\JL$ is a complemented
copy of~$\ell_2$.
\end{Cor}

\begin{Proof}
Let $M$ be an \infdim{} reflexive subspace of~$\JL$.
Then $M \cap V$ is finite-dimensional
and, by Lemma~\ref{sub-3sp-lemma}, $M + V$ is closed.
Passing to a finite-codimensional subspace~$N$ of~$M$, the restriction
$Q_V J_N$ is an isomorphism. Since $Q_V(N)$ is complemented
in~$\ell_2(\Gamma)$, $N$ is complemented in~$\JL$ by Lemma~\ref{useful}(i),
and then so is~$M$.
\end{Proof}

Using subprojective and superprojective operators, it is possible to prove
that the space~$\JL$ is both subprojective and superprojective.

\begin{Lemma}\label{disjoint}
Let $\{\, \gamma_i : 1\leq i\leq n \,\}$ be a finite subset of~$\Gamma$ and
let $M_k = N_{\gamma_k} \setminus \bigcup_{i=1}^{k-1} N_{\gamma_i}$ for every
$1\leq k \leq n$. Then $(\chi_{M_k})_{i=1}^n$ is equivalent to the unit
vector basis of~$\ell_2^n$.
\end{Lemma}

\begin{Proof}
For every $1\leq k\leq n$, define $F_k = \bigcup_{i=1}^{k-1}
(N_{\gamma_k} \cap N_{\gamma_i})$, which is a finite set, and note that
$N_{\gamma_k}$ is the disjoint union of $M_k$ and~$F_k$,
so $\chi_{M_k} = \chi_{N_{\gamma_k}} - \chi_{F_k}$, which is well defined
in~$\JL$. Since the sets $(M_k)_{k=1}^n$ are pairwise disjoint,
$\bigl\| \sum_{k=1}^n a_k \chi_{M_k} \bigr\|_\infty =
\|(a_k)_{k=1}^n\|_\infty$ and
$\bigl\| \sum_{k=1}^n a_k \chi_{M_k} \bigr\|_\JL =
\|(a_k)_{k=1}^n\|_2$ for every $(a_k)_{k=1}^n$.
\end{Proof}

As a consequence, the same applies to any countable subset of~$\Gamma$
with respect to~$\ell_2$.

We will need the following result, which was essentially proved in
\cite[Theorem~2.2]{diaz-fernandez}. We include a proof here because our
statement is different and also encompasses complex spaces.

\begin{Prop}\label{c0-no-l1}
Let $X$ be a Banach space that does not contain any copies of~$\ell_1$.
Then every copy of~$c_0$ in~$X$ contains another copy of~$c_0$ complemented
in~$X$.
\end{Prop}

\begin{Proof}
Let $(x_n)_{n\in\N}$ be a sequence in~$X$ equivalent to the unit vector basis
of~$c_0$ and take a bounded sequence $(x^*_n)_{n\in\N}$ in~$X^*$ such that
$\<x^*_i,x_j\> = \delta_{ij}$ for every $i$, $j\in\N$. It is easy to check
that $(x^*_n)_{n\in\N}$ is equivalent to the unit vector basis of~$\ell_1$.

Since $X$ does not contain any copies of~$\ell_1$, there exists a normalised
weak$^*$ null block sequence of~$(x^*_n)_{n\in\N}$
(\cite[Theorem~1(a)]{hagler-johnson} for the real case,
\cite[Appendix~A]{acuaviva} for the complex case).
Write $y^*_k = \sum_{i=m_k}^{m_{k+1}-1} a_i x^*_i$ for such a sequence
and let $\eps_i = a_i/|a_i|$, or $\eps_i=1$ if $a_i=0$, for every $i\in\N$.
Then $y_k = \sum_{i=m_k}^{m_{k+1}-1} \eps_i x_i$ defines a sequence
$(y_k)_{k\in\N}$ in $[x_n:n\in\N]$ equivalent to the unit vector basis
of~$c_0$ such that $P(x) = \sum_{k=1}^\infty \<y^*_k,x\> y_k$ is a
projection in~$X$ onto~$[y_k:k\in\N]$.
\end{Proof}

\begin{Prop}
\leavevmode
\begin{enumerate}
\itemn1 $J_V$ is subprojective.
\itemn2 $Q_V$ is superprojective.
\itemn3 $\JL$ is both subprojective and superprojective.
\itemn4 $\JL^*$ is subprojective but not superprojective.
\itemn5 $\JL^{**}$ is neither subprojective nor superprojective.
\end{enumerate}
\end{Prop}

\begin{Proof}
(i)~$\JL$ does not contain any copies of~$\ell_1$, as neither $c_0$
nor~$\ell_2(\Gamma)$ contain copies of~$\ell_1$, and not containing any copies
of~$\ell_1$ is a three-space property \cite[Theorem~3.2.d]{3-space}. Then
every closed \infdim{} subspace of~$V\simeq c_0$ contains another
copy of~$c_0$ that is complemented in~$\JL$ by Proposition~\ref{c0-no-l1}
and $J_V$ is subprojective by Proposition~\ref{wrt}.

(ii)~Let $M$ be a closed \infcodim{} subspace of~$\JL$
containing~$V$. Then $\JL/M$ is a quotient of $\JL/V \equiv
\ell_2(\Gamma)$, so, taking a bigger~$M$
if necessary, we can assume that $\JL/M$ is separable. Since $Q_V(M)$
is closed, we can consider the decomposition $\ell_2(\Gamma) =
Q_V(M) \oplus Q_V(M)^\bot$.

Let $\{\, e_\gamma : \gamma\in\Gamma \,\}$ be the basis of~$\ell_2(\Gamma)$.
Since $Q_V(M)^\bot$ is separable, there exists a sequence of different
points $(\gamma_n)_{n\in\N}$ in~$\Gamma$ such that $Q_V(M)^\bot \subseteq
[e_{\gamma_n}:n\in\N]$. Consider an orthonormal basis of~$Q_V(M)^\bot$;
then that basis is a normalised weakly null sequence, so it has a subsequence
$(f_k)_{k\in\N}$ that is equivalent to a block basis
of~$(e_{\gamma_n})_{n\in\N}$, i.e., there exist $1 = m_1 < m_2 < m_3 < \cdots$
in~$\N$ and a sequence of scalars $(b_n)_{n\in\N}$ such that $y_k =
\sum_{i=m_k}^{m_{k+1}-1} b_i e_{\gamma_i}$ satisfies $\|y_k - f_k\| < 2^{-2k}$
for every $k\in\N$.

Now, for every $n\in\N$, define $F_n = \bigcup_{i=1}^{n-1}
(N_{\gamma_n} \cap N_{\gamma_i})$ and $M_n = N_{\gamma_n} \setminus F_n$,
so that $(\chi_{M_n})_{n\in\N}$ is equivalent to the unit vector basis
of~$\ell_2$ by Lemma~\ref{disjoint} and $Q_V(\chi_{M_n}) =
Q_V(\chi_{N_{\gamma_n}}) = e_{\gamma_n}$ for every $n\in\N$,
which means that the restriction of~$Q_V$ to $[\chi_{M_k}:k\in\N]$ is
an isomorphism onto $[e_{\gamma_n}:n\in\N]$.
Define $x_k = \sum_{i=m_k}^{m_{k+1}-1} b_i \chi_{M_i} \in \JL$ for
every $k\in\N$, so that $Q_V(x_k) = y_k$ for every $k\in\N$. Then
the restriction of~$Q_V$ to~$[x_k:k\in\N]$ is an isomorphism
onto~$[y_k:k\in\N]$ and $\ell_2(\Gamma) = [y_k:k\in\N] \oplus N$ with $N$
a closed subspace containing $Q_V(M)$, and $\JL = [x_k:k\in\N] \oplus
Q_V^{-1}(N)$. Hence $Q_V^{-1}(N)$ is a complemented \infcodim{}
subspace containing~$M$.

(iii) $Q_V$ is subprojective and $J_V$ is superprojective
by Corollary~\ref{space-operator}, as $\JL/V \equiv \ell_2(\Gamma)$ is
subprojective and $V = c_0$ is superprojective, so $\JL$ is both
subprojective and superprojective by (i), (ii) and Theorem~\ref{3sp}.

(iv) This follows from $\JL^* \simeq \ell_1 \oplus \ell_2(\Gamma)$
\cite[Proposition~2.2]{oikhberg-spinu}
\cite[Proposition~4.1]{superprojective}, as
$\ell_1$ and $\ell_2(\Gamma)$ are subprojective, but $\ell_1$
is not superprojective.

(v) $\JL^{**}$ contains a (complemented) copy of~$\ell_\infty$.
\end{Proof}

There is a second example in~\cite{wcg}, which is the closed subspace~$X_\JL$
generated by $V \cup \{\, \phi_\gamma : \gamma\in\Gamma \,\} \cup \{\chi_\N\}$
in~$\ell_\infty$. Since $X_\JL$ is a commutative Banach algebra with
the pointwise multiplication, it is isometric to some $C(K)$ space by
Gelfand's representation theorem \cite[Example~2]{wcg}.
It is easy to check that $X_\JL/V$ is isometric to $c_0(\Gamma)$, which
is WCG because the natural inclusion $\ell_2(\Gamma) \maparrow c_0(\Gamma)$
has dense range. However, since weakly compact subsets of~$\ell_\infty$
are separable, $X_\JL$ is not WCG, so $V$ is not complemented in~$X_\JL$.

Recall that a Banach space~$X$ is said to have property~$(V)$
if every non-weakly compact operator $\map TXY$ is an isomorphism on a
subspace of~$X$ isomorphic to~$c_0$. It is well known that $C(K)$ spaces
have property~$(V)$, so $X_\JL$ has property~$(V)$.

\begin{Prop}
$X_\JL$ is subprojective and superprojective.
\end{Prop}

\begin{Proof}
$X_\JL$ is hereditarily~$c_0$ because being hereditarily~$c_0$
is a three-space property \cite[Theorem~3.2.e]{3-space}. As such,
any closed \infdim{} subspace of~$V = c_0$ contains a further
subspace isomorphic to~$c_0$ and complemented in~$X_\JL$ by
Proposition~\ref{c0-no-l1}, which means that $J_V$ is subprojective.
As in the case of~$\JL$, $Q_V$ is subprojective by
Corollary~\ref{space-operator}, as $X_\JL/V \equiv c_0(\Gamma)$
is subprojective, so $X_\JL$ is subprojective by Theorem~\ref{3sp}.

For the superprojectivity, we will prove that $Q_V$ is superprojective.
First of all, $\ell_1(\Gamma)$ does not have any reflexive subspaces,
so $X_\JL/V \equiv c_0(\Gamma)$ does not have any reflexive quotients.
Let $W$ be a closed \infcodim{} subspace of~$X_\JL$ containing~$V$;
then $Q_W$ is not weakly compact, and $X_\JL$ has property~$(V)$, so there
exists a subspace~$W$ of~$X_\JL$ isomorphic to~$c_0$ such that $Q_W$ is an
isomorphism on~$M$. Now again $X_\JL/W$ does not have any reflexive
quotients, so it does not contain~$\ell_1$ and we can assume that $Q_W(M)$ is
complemented by Proposition~\ref{c0-no-l1}. Write $X_\JL/W = Q_W(M) \oplus N$;
then $X_\JL = M \oplus Q_W^{-1}(N)$, where $Q_W^{-1}(N)$ contains~$W$.

Finally, again $J_V$ is superprojective by Corollary~\ref{space-operator},
so $X_\JL$ is superprojective by Theorem~\ref{3sp}.
\end{Proof}

\end{section}


\begin{thebibliography}{99}

\bibitem{acuaviva}
A.\ Acuaviva.
\emph{Factorizations and minimality of the Calkin algebra norm for
  $C(K)$-spaces.}
arXiv:2408.11132v1, 2024.

\bibitem{improjective}
P.\ Aiena, M.\ Gonz\'alez.
\emph{On inessential and improjective operators.}
Studia Math.\ 131 (1998) 271--287.

\bibitem{kalton-albiac}
F.\ Albiac, N.\ J.\ Kalton.
\emph{Topics in Banach Space Theory.}
Graduate Texts in Mathematics, 233.
Springer, 2006.

\bibitem{argyros-raikoftsalis}
S.\ A.\ Argyros, T.\ Raikoftsalis.
\emph{The cofinal problem property of the reflexive indecomposable Banach
  spaces.}
Ann.\ Inst.\ Fourier (Grenoble) 62 (2012), 1--45.

\bibitem{jl-weak-star}
A.\ Avil\'es, G.\ Mart\'\i nez-Cervantes, J.\ Rodr\'\i guez.
\emph{Weak*-sequential properties of Johnson-Lindenstrauss spaces.}
J.\ Funct.\ Anal.\ 276 (10) (2019), 3051--3066.

\bibitem{Jsum}
S.\ F.\ Bellenot.
\emph{The $J$-sum of Banach spaces.}
J.\ Funct.\ Anal.\ 48 (1) (1982) 95--106.

\bibitem{3-space}
J.\ M.~F.\ Castillo, M.\ Gonz\'alez.
\emph{Three-space problems in Banach space theory.}
Lecture Notes in Mathematics, 1667.
Springer-Verlag, 1997.

\bibitem{diaz-fernandez}
S.\ D\'\i az and A.\ Fern\'andez.
\emph{Reflexivity in Banach lattices.}
Arch.\ Math.\ 63 (6) (1994), 549--552.

\bibitem{perturbation}
M.\ Gonz\'alez.
\emph{The perturbation classes problem in Fredholm theory.}
J.\ Funct.\ Anal.\ 200 (2003) 65--70.

\bibitem{tauberian}
M.\ Gonz\'alez, A.\ Mart\'\i nez-Abej\'on.
\emph{Tauberian Operators.}
Operator Theory: Advances and Applications, 194.
Birk\"auser, 2010.

\bibitem{perturbation-sub-sup}
M.\ Gonz\'alez, A.\ Mart\'\i nez-Abej\'on, M.\ Salas-Brown.
\emph{Perturbation classes for semi-Fredholm operators on subprojective
  and superprojective Banach spaces.}
Ann.\ Acad.\ Sci.\ Fenn.\ Math.\ 36 (2011) 481--491.

\bibitem{j-sum-subprojective}
M.\ Gonz\'alez, J.\ Pello.
\emph{Projections in the $J$-sums of Banach spaces.}
Rev.\ Real Acad.\ Cienc.\ Exactas F\'\i s.\ Nat.\ Ser.\ A-Mat.\ 118 (4) (2024), 162.

\bibitem{superprojective}
M.\ Gonz\'alez, J.\ Pello.
\emph{Superprojective Banach spaces.}
J.\ Math.\ Anal.\ Appl.\ 437 (2) (2016) 1140--1151.

\bibitem{G-P-Salas}
M.\ Gonz\'alez, J.\ Pello, M.\ Salas-Brown.
\emph{The perturbation classes problem for subprojective and superprojective
  Banach spaces.}
J.\ Math.\ Anal.\ Appl.\ 489 (2) (2020) 124191, 5~pp.

\bibitem{hagler-johnson}
J.\ Hagler, W.\ B.\ Johnson.
\emph{On Banach spaces whose dual balls are not weak* sequentially compact.}
Israel J.\ Math.\ 28 (4) (1977), 325--330.

\bibitem{wcg}
W.\ B.\ Johnson, J.\ Lindenstrauss.
\emph{Some remarks on weakly compactly generated Banach spaces.}
Israel J.\ Math.\ 17 (2) (1974) 219--230.

\bibitem{kalton-peck}
N.\ J.\ Kalton, N.\ T.\ Peck.
\emph{Twisted sums of sequence spaces and the three space problem.}
Trans.\ Amer.\ Math.\ Soc.\ 255 (1979) 1--30.

\bibitem{kato}
T.\ Kato.
\emph{Perturbation Theory for Linear Operators.}
Springer, 1966.

\bibitem{oikhberg-spinu}
T.\ Oikhberg, E.\ Spinu.
\emph{Subprojective Banach spaces.}
J.\ Math.\ Anal.\ Appl.\ 424 (2015) 613--635.

\bibitem{whitley}
R.\ J.\ Whitley.
\emph{Strictly singular operators and their conjugates.}
Trans.\ Amer.\ Math.\ Soc.\ 113 (1964) 252--261.

\end{thebibliography}
\end{document}